# ON THE LOCATION AND CLASSIFICATION OF ALL PRIME NUMBERS

Leopoldo Garavaglia, Aranjuez, ESPAÑA, and
Mario Garavaglia, Universidad Nacional de La Plata and Centro de Investigaciones Opticas, La Plata, ARGENTINA.

**Abstract:** We will describe an algorithm to arrange all the positive and negative integer numbers. This array of numbers permits grouping them in six different Classes, α, β, γ, δ, ε, and ζ. Particularly, numbers belong to Class α are defined as α = 1 + 6 n, and those of Class β, as β = 5 + 6 n, where n = 0, ± 1, ± 2, ± 3, ± 4, … These two Classes α and β, contain: i) all prime numbers, except + **2**, –**2** and ± **3**, which belong to ε, δ, and γ Classes, respectively, and ii) all the other odd numbers, except those that are multiple of ± **3**, according to the sequence ± **9**, ± **15**, ± **21**, ± **27**, … Besides, products between numbers of the Class α, and also those between numbers of the Class β, generates numbers belonging to the Class α. On the other side, products between numbers of Class α with numbers of Class β, result in numbers of Class β. Then, both Classes α and β include: i) all the prime numbers except ± **2** and ± **3**, and ii) all the products between α numbers, as **α. α'**; all the products between β numbers, as **β.β'**; and also all the products between numbers of Classes α and β, as **α.β**, which necessarily are composite numbers, whose factorization is completely determined.

**Key words:** Prime numbers. Arrange algorithm. Matrix of products. 2D representations. Location and classification of prime numbers. Factorization.

## 1. INTRODUCTION [1]

It is well recognized that the idea of a decimal algorithm to methodical orderly array of integer numbers was introduced in India in the VI Century BC; such arrangement incorporates also the zero. However, to our purpose we will use another orderly array of even integers numbers –or natural numbers– as follow:

$$
\begin{array}{ccccccccc}
1 & 2 & 3 & 4 & 5 & 6 & 7 & 8 & 9 \\
10 & 11 & 12 & 13 & 14 & 15 & 16 & 17 & 18 \\
19 & 20 & 21 & 22 & 23 & 24 & 25 & 26 & 27 \ldots
\end{array}
$$

It is easy to prove that all numbers in the first column belong to the Digital Root, that is:

$$DR(\mathbf{1}) = \mathbf{1}; \ DR(\mathbf{10}) \equiv \mathbf{1 + 0 = 1}; \ DR(\mathbf{19}) \equiv \mathbf{1 + 9 = 1}; \ldots$$

Besides, all the numbers in the second column belong to the DR(**2**) = **2**, and so on and so forth. Then, we will designate all numbers in the first column as numbers of Type **1** or **a**, those of the second column as belonging to the Type **2** or **b** numbers, and so on. In the Table I is represented our arrangement, which we named First Array (A1).

| **a** | **b** | **c** | **d** | **e** | **f** | **g** | **h** | **i** |
|---|---|---|---|---|---|---|---|---|
| **1** | **2** | **3** | **4** | **5** | **6** | **7** | **8** | **9** |
| **10** | **11** | **12** | **13** | **14** | **15** | **16** | **17** | **18** |
| **19** | **20** | **21** | **22** | **23** | **24** | **25** | **26** | **27 . . .** |

Table I: Partial representation of A1.

The contained information in A1 permits to elaborate a general matrix to represent the results of the products of all the Type numbers in Table I.

## 2. MATRIX OF PRODUCTS M [1]

Table II describes the algorithm of the Matrix of Products M.

| M | **X** | a | b | c | d | e | f | g | h | i |
|---|---|---|---|---|---|---|---|---|---|---|
| | a | a | b | c | d | e | f | g | h | i |
| | b | b | d | f | h | a | c | e | g | i |
| | c | c | f | i | c | f | i | c | f | i |
| | d | d | h | c | g | b | f | a | e | i |
| | e | e | a | f | b | g | c | h | d | i |
| | f | f | c | i | f | c | i | f | c | i |
| | g | g | e | c | a | h | f | d | b | i |
| | h | h | g | f | e | d | c | b | a | i |
| | i | i | i | i | i | i | i | i | i | i |

Table II: The Matrix of Products M.

## 3. LOCATION OF PRIME NUMBERS IN A1 [1]

It is easy to observe that all prime numbers –except **3**– are located in A1 as numbers of Types **a**, **b**, **d**, **e**, **g**, or **h**. The only even prime number 2 belongs to the second column which corresponds to the Type **b** numbers.

According to the theorem: "the only divisors of 1 are ± **1**", it results that prime numbers **p** different to **0** and ± **1**, that can only be divided by ± **1** and ± **p**. By using this definition the previous array A1 can be extended to include the **0** and all the negative integers' numbers, as follows:

$$0 \ -1 \ -2 \ -3 \ -4 \ -5 \ -6 \ -7 \ -8$$
$$-9 \ -10 \ -11 \ -12 \ -13 \ -14 \ -15 \ -16 \ -17$$
$$-18 \ -19 \ -20 \ -21 \ -22 \ -23 \ -24 \ -25 \ -26 \ldots$$

We will refer to this array as the Second Array (A2). Table III shows three lines of the positive integer numbers and three lines of the negative integer numbers, just crossing the **0**; the **0** is a Type **i** number.

|   | a | b | c | d | e | f | g | h | i |
|---|---|---|---|---|---|---|---|---|---|
| ... | -26 | -25 | -24 | -23 | -22 | -21 | -20 | -19 | -18 |
|   | -17 | -16 | -15 | -14 | -13 | -12 | -11 | -10 | -9 |
|   | -8 | -7 | -6 | -5 | -4 | -3 | -2 | -1 | 0 |
|   | 1 | 2 | 3 | 4 | 5 | 6 | 7 | 8 | 9 |
|   | 10 | 11 | 12 | 13 | 14 | 15 | 16 | 17 | 18 |
|   | 19 | 20 | 21 | 22 | 23 | 24 | 25 | 26 | 27 ... |

Table III: The transition from A1 to A2 through the **0**.

## 4. THE TRANSITION AND THE MATRIX OF PRODUCTS [1]

Note that crossing the transition from A1 to A2 the column of the Type **h** numbers contains the sequence …, 26, 17, 8, −1, −10, −19, …, which at first sight looks unacceptable. However, the Matrix of Products demonstrates that such a sequence is necessarily logic. In fact, the multiplication 8 x (−1), which are both **h** numbers, is equal to −8, which is an **a** number, just as in the case of the multiplication 8 x 8 = 64, in accordance to the Matrix M of Products.

In A2 the **0** occupy the column of numbers Type **i**, while the negative prime numbers, except −3 which is an **f** number, belong to the Type **h**, **g**, **e**, **d**, **b**, or **a** numbers; −2, the only negative even prime numbers is of the Type **g**.

## 5. FROM A2 TO A NEW ARRAY A3

We rearranged A2 in order to distribute the nine columns of numbers **a**, **b**, **c**, **d**, **e**, **f**, **g**, **h**, and **i**, in the following three columns:

|   |   |   |
|---|---|---|
| **a** | **b** | **c** |
| **d** | **e** | **f** |
| **g** | **h** | **i** ; |

we named this array as A3. Table IV contains the new arrangement A3.

| | | |
|---|---|---|
| . . . **a** -26 | **b** -25 | **c** -24 |
| **d** -23 | **e** -22 | **f** -21 |
| **g** -20 | **h** -19 | **i** -18 |
| **a** -17 | **b** -16 | **c** -15 |
| **d** -14 | **e** -13 | **f** -12 |
| **g** -11 | **h** -10 | **i** - 9 |
| **a** - 8 | **b** - 7 | **c** - 6 |
| **d** - 5 | **e** - 4 | **f** - 3 |
| **g** - 2 | **h** - 1 | **i** 0 |
| **a** 1 | **b** 2 | **c** 3 |
| **d** 4 | **e** 5 | **f** 6 |
| **g** 7 | **h** 8 | **i** 9 |
| **a** 10 | **b** 11 | **c** 12 |
| **d** 13 | **e** 14 | **f** 15 |

|   |   |   |   |   |   |
|---|---|---|---|---|---|
| g | 16 | h | 17 | i | 18 |
| a | 19 | b | 20 | c | 21 |
| d | 22 | e | 23 | f | 24 |
| g | 25 | h | 26 | i | 27 ... |

Table IV. Array A3.

As A2, the A3 array contains all the positive and negative integer numbers. So, from A3 we will now perform two new arrangements of numbers: OA for odd numbers and EA for even numbers, as follows:

|   |   |   |   |   |   |
|---|---|---|---|---|---|
| **Odd** | g | 7 | b | 11 | i | 9 |
| **Even** | a | 10 | h | 8 | c | 12. |

Table V includes both arrays OA and y OE.

|   | **Odd** |   |   |   | **Even** |   |
|---|---|---|---|---|---|---|
| ... g -29 | b -25 | i -27 | ... a -26 | h -28 | c -24 |
| d -23 | h -19 | f -21 | g -20 | e -22 | i -18 |
| a -17 | e -13 | c -15 | d -14 | b -16 | f -12 |
| g -11 | b - 7 | i - 9 | a - 8 | h -10 | c - 6 |
| d - 5 | h - 1 | f - 3 | g - 2 | e - 4 | i   0 |
| a   1 | e   5 | c   3 | d   4 | b   2 | f   6 |
| g   7 | b  11 | i   9 | a  10 | h   8 | c  12 |
| d  13 | h  17 | f  15 | g  16 | e  14 | i  18 |
| a  19 | e  23 | c  21 | d  22 | b  20 | f  24 |
| g  25 | b  29 | i  27 ... | a  28 | h  26 | c  30 ... |

Table V. The arrays OA and EA of numbers.

Note that in the three columns of OA as well as in the three columns of EA the difference between numbers are always six units. Then, we will separate those numbers in six Classes, as we show in Table VI.

| **Classes:** | α | β | γ | δ | ε | ζ |
|---|---|---|---|---|---|---|
| | d  - 5 | h  - 1 | f  - 3 | g  - 2 | e  - 4 | i  0 |
| | a  1 | e  5 | c  3 | d  4 | b  2 | f  6 |

Table VI. Each column in Table V are represented by different Classes of Numbers **α**, **β**, **γ**, **δ**, **ε**, and **ζ**.

Numbers in each Classes **α**, **β**, **γ**, **δ**, **ε**, and **ζ**, can be expressed by the following equations:

$\alpha = 1 + 6n$,   [1]     $\beta = 5 + 6n$,   [2]     $\gamma = 3 + 6n$,   [3]

$\delta = 4 + 6n$,   [4]     $\varepsilon = 2 + 6n$,   [5]     $\zeta = 6 + 6n$,   [6]

being **n** = 0, ± 1, ± 2, ± 3, ± 4, ...

### 6. LOCATION OF ALL PRIME NUMBERS

According to Table VI positive and negative integer numbers in each different Classes **α**, **β**, **γ**, **δ**, **ε**, and **ζ**, belongs only to certain Types of numbers, as follows:

- Numbers **α**, can only be of Types **a**, **d**, **g**,
- Numbers **β**, can only be of Types **b**, **e**, **h**,
- Numbers **γ**, can only be of Types **c**, **f**, **i**,
- Numbers **δ**, can only be of Types **a**, **d**, **g**,
- Numbers **ε**, can only be of Types **b**, **e**, **h**,
- Numbers **ζ**, can only be of Types **c**, **f**, **i**.

Then, the only positive even prime number **2** is an **ε b** number, while **– 2**, the only negative even prime number, is a **δ g** number. Besides, **+ 3** is a **γ c**, the only positive odd prime number, while **– 3** is also a **γ** number, but of Type **f**, the only negative even prime number. All the others **c**, **f**, and **i** numbers are multiple numbers of ± **3**, so they are composite numbers. Then, we must conclude that all the prime numbers –except ± **2** and ± **3**– are located in Classes **α** and **β**. But, as in both Classes **α** y **β** are located others positive and negative numbers which are not prime numbers, we must also conclude that they are necessarily composite numbers.

This last conclussion can be demonstrated by following two ways: i) by using the Matrix M of Products, or ii) by multipling numbers of different Classes **α**, **β**, **γ**, **δ**, **ε**, and **ζ**.

**6.1.** Demonstration by using the Matrix M of Products:

From the matrix M in Table II we will extract the terms **a**, **d**, and **g**, and also the terms **b**, **e**, and **h**, to build up the three following matrices $M_1$, $M_2$, and $M_3$:

| $M_1$ | X | a | d | g |
|---|---|---|---|---|
| | a | a | d | g |
| | d | d | g | a |
| | g | g | a | d, |

| $M_2$ | X | b | e | h |
|---|---|---|---|---|
| | b | d | a | g |
| | e | a | g | d |
| | h | g | d | a, |

| $M_3$ | X | a | d | g |
|---|---|---|---|---|
| | b | b | h | e |
| | e | e | b | h |
| | h | h | e | b. |

From one side, products in $M_1$ involving numbers of the Type **a**, **d**, and **g**, result in numbers of the same Type **a**, **d**, and **g**. Also, the products in $M_2$ involving numbers of the Type **b**, **e**, and **h**, result in numbers of the same Type **b**, **e**, y **h**. On the other side, products in $M_3$ involving numbers of the Type **a**, **d**, and **g**, and numbers of the Type **b**, **e**, y **h**, result in numbers of the Type **b**, **e**, y **h**. Then, we must conclude that all products represented inside matrices $M_1$, $M_2$, and $M_3$ belong to both Classes **α** y **β**, and they are necessarily composite numbers.

Going back to the matrix M in Table II, we will find that all the other products involving numbers of the Type **a**, **b**, **d**, **e**, **g**, and **h**, with numbers of the Type **c**, **f**, and **i**, always result in numbers of the Type *c*, *f*, and *i*, as it is shown in *italic* in Table VII.

| M | X | a | b | c | d | e | f | g | h | i |
|---|---|---|---|---|---|---|---|---|---|---|
|   | **a** | a | b | *c* | d | e | *f* | g | h | *i* |
|   | **b** | b | d | *f* | h | a | *c* | e | g | *i* |
|   | **c** | *c* | *f* | *i* | *c* | *f* | *i* | *c* | *f* | *i* |
|   | **d** | d | h | *c* | g | b | *f* | a | e | *i* |
|   | **e** | e | a | *f* | b | g | *c* | h | d | *i* |
|   | **f** | *f* | *c* | *i* | *f* | *c* | *i* | *f* | *c* | *i* |
|   | **g** | g | e | *c* | a | h | *f* | d | b | *i* |
|   | **h** | h | g | *f* | e | d | *c* | b | a | *i* |
|   | **i** | *i* | *i* | *i* | *i* | *i* | *i* | *i* | *i* | *i* |

Table VII. Matrix M of Products in which we emphasized in *italic* results of products involving numbers of the Type **a**, **b**, **d**, **e**, **g**, and **h**, with numbers of the Type **c**, **f**, and **i**.

Then, the application of the matrix M confirms that the products involving numbers of the Type **a**, **b**, **d**, **e**, **g**, and **h**, with numbers of the Type **c**, **f**, and **i**, result in numbers of the Type *c*, *f*, and *i*, in *italic* in Table VII, which are always composite numbers whose factorization is completely determined.

**6.2.** Demonstration by using equations [1] to [6]:

The application of the equations [1] to [6] to perform products between numbers from the different Classes α, β, γ, δ, ε, y ζ, is the other way we propose to demonstrate that all prime numbers –except ± **2** and ± **3**– are located in Classes **α** and **β**, and their products **α.α'**, **β.β'**, and **α.β** are composite numbers. Besides, we will demonstrate that products between numbers from the different Classes α, β, γ, δ, ε, y ζ, as **α. γ, α. δ, α. ε, α. ζ, β. γ, β. δ, β. ε, β. ζ, γ. δ, γ. ε, γ. ζ, δ.ε, δ.ζ**, and **ζ.ζ**, do not result in numbers of Classes **α** or **β**.

Let be, $\alpha = 1 + 6 n_\alpha$, and $\alpha' = 1 + 6 n'_\alpha$, then, the product:

$\alpha.\alpha' = (1 + 6 n_\alpha).(1 + 6 n'_\alpha) = 1 + 6 n_\alpha + 6 n'_\alpha + 36 n_\alpha n'_\alpha = 1 + 6 (n_\alpha + n'_\alpha + 6 n_\alpha n'_\alpha)$

$= \alpha'' = 1 + 6 n''_\alpha,$

so, the product **α.α'** = **α''**, which belongs to the Class **α**, being $n''_\alpha = n_\alpha + n'_\alpha + 6 n_\alpha n'_\alpha$.

Besides, let be, $\beta = 5 + 6 n_\beta$, and $\beta' = 5 + 6 n'_\beta$, then the product:

$\beta.\beta' = (5 + 6 n_\beta).(5 + 6 n'_\beta) = 25 + 30 n_\beta + 30 n'_\beta + 36 n_\beta n'_\beta =$

$1 + 24 + 30 n_\beta + 30 n'_\beta + 36 n_\beta n'_\beta = 1 + 6 (4 + 5 n_\beta + 5 n'_\beta + 6 n_\beta n'_\beta) =$

$\alpha = 1 + 6 n''_\alpha,$

so, the product **β.β'** = **α**, being $n''_\alpha = 4 + 5 n_\beta + 5 n'_\beta + 6 n_\beta n'_\beta$.

In the case that the product combines a number belongs to Class **α** and a number from the **β** Class, defined as $\alpha = 1 + 6 n_\alpha$, and $\beta = 5 + 6 n_\beta$, respectively, it results:

$\alpha.\beta = (1 + 6 n_\alpha).(5 + 6 n_\beta) = 5 + 30 n_\alpha + 6 n_\beta + 36 n_\alpha n_\beta = 5 + 6 (5 n_\alpha + n_\beta + 6 n_\alpha n_\beta) =$

$\beta' = 5 + 6 n'_\beta,$

so, the product **α.β** = **β'**, being $n'_\beta = 5 n_\alpha + n_\beta + 6 n_\alpha n_\beta$.

Using the same procedure to perform the products **α. γ, α. δ, α. ε, α. ζ, β. γ, β. δ, β. ε, β. ζ, γ. δ, γ. ε, γ. ζ, δ.ε, δ.ζ**, and **ζ.ζ**, we will prove that they do not result in numbers of Classes **α** or **β**.

In fact:

- $\alpha. \gamma = (1 + 6 n_\alpha).( 3 + 6 n_\gamma) = \gamma' = 3 + 6 n'_\gamma,$

being $n'_\gamma = n_\gamma + 3 n_\alpha + 6 n_\alpha n_\gamma$.

- $\alpha. \delta = \delta' = (1 + 6 n_\alpha).( 4 + 6 n_\delta) = \delta' = 4 + 6 n'_\delta,$

being $n'_\delta = n_\delta + 4 n_\alpha + 6 n_\alpha n_\delta$.

- α. ε = (1 + 6 n$_α$).( 2 + 6 n$_ε$) = ε' = 2 + 6 n'$_ε$,

being n'$_ε$ = n$_ε$ + 2 n$_α$ + 6 n$_α$ n$_ε$.

- α. ζ = (1 + 6 n$_α$).( 6 + 6 n$_ζ$) = ζ' = 6 + 6 n'$_ζ$,

being n'$_ζ$ = n$_ζ$ + 6 n$_α$ + 6 n$_α$ n$_ζ$.

- β. γ = (5 + 6 n$_β$).( 3 + 6 n$_γ$) = γ' = 3 + 6 n'$_γ$,

being n'$_γ$ = 2 + 3 n$_β$ + 5 n$_γ$ + 6 n$_β$ n$_γ$.

- β. δ = (5 + 6 n$_β$).( 4 + 6 n$_δ$) = ε = 2 + 6 n$_ε$,

being n$_ε$ = 3 + 4 n$_β$ + 5 n$_δ$ + 6 n$_β$ n$_δ$.

- β. ε = (5 + 6 n$_β$).( 2 + 6 n$_ε$) = δ = 4 + 6 n$_δ$,

being n$_δ$ = 1 + 2 n$_β$ + 5 n$_ε$ + 6 n$_β$ n$_ε$.

- β. ζ = (5 + 6 n$_β$).( 6 + 6 n$_ζ$) = ζ' = 6 + 6 n'$_ζ$,

being n'$_ζ$ = 4 + 6 n$_β$ + 5 n$_ζ$ + 6 n$_β$ n$_ζ$.

- γ. δ = (3 + 6 n$_γ$).( 4 + 6 n$_δ$) = ζ = 6 + 6 n$_ζ$,

being n$_ζ$ = 1 + 3 n$_δ$ + 4 n$_γ$ + 6 n$_δ$ n$_γ$.

- γ. ε = (3 + 6 n$_γ$).( 2 + 6 n$_ε$) = ζ = 6 + 6 n$_ζ$,

being n$_ζ$ = 2 n$_γ$ + 3 n$_ε$ + 6 n$_γ$ n$_ε$.

- γ. ζ = (3 + 6 n$_γ$).( 6 + 6 n$_ζ$) = ζ' = 6 + 6 n'$_ζ$,

being n$_ζ$ = 2 + 6 n$_γ$ + 3 n$_ζ$ + 6 n$_γ$ n$_ζ$.

- δ . ε = ( 4 + 6 n$_δ$).(2 + 6 n$_ε$) = ε' = 2 + 6 n'$_ε$,

being n'$_ε$ = 1 + 2 n$_δ$ + 4 n$_ε$ + 6 n$_δ$ n$_ε$.

- δ . ζ = ( 4 + 6 n$_δ$).(6 + 6 n$_ζ$) = ζ' = 6 + 6 n'$_ζ$,

being n'$_ζ$ = 3 + 6 n$_δ$ + 4 n$_ζ$ + 6 n$_δ$ n$_ζ$.

- ζ . ζ' = ( 6 + 6 n$_ζ$).(6 + 6 n'$_ζ$) = ζ'' = 6 + 6 n''$_ζ$,

being n''$_ζ$ = 5 + 6 n$_ζ$ + 6 n'$_ζ$ + 6 n$_ζ$ n'$_ζ$.

Then, from demonstrations in **6.1.** and in **6.2.**, we conclude that all the products between numbers from the Class **α** and numbers from the Class **β**, always result in numbers that belong to Class **α**, while all the products involving numbers from Class **α** and Class **β**, always result in numbers that belong to Class **β**, which confirms that: i) all prime numbers – except **± 2** y **± 3**– are located inside Classes **α** and **β**, and ii) all the others positive and negative numbers which are not prime numbers which are located in Classes **α** y **β**, are necessarily composite numbers, whose factorization is completely determined.

**7. LOCATION OF PRIME NUMBERS AND COMPOSITE NUMBERS UP TO 605**

According to the results obtained in Section **6.** the next Table VIII includes all **α** and **β** numbers up to **605**; when one of them is a composite number its expression as a product of prime numbers is given.

| n | Type | Class α | Type | Class β |
|---|---|---|---|---|
| 0 | a | 1 | e | 5 |
| 1 | g | 7 | b | 11 |
| 2 | d | 13 | h | 17 |
| 3 | a | 19 | e | 23 |
| 4 | g | 25 = 5 . 5 | b | 29 |
| 5 | d | 31 | h | 35 = 5 . 7 |
| 6 | a | 37 | e | 41 |
| 7 | g | 43 | b | 47 |
| 8 | d | 49 = 7 . 7 | h | 53 |
| 9 | a | 55 = 5 . 11 | e | 59 |
| 10 | g | 61 | b | 65 = 5 . 13 |
| 11 | d | 67 | h | 71 |
| 12 | a | 73 | e | 77 = 7 . 11 |
| 13 | g | 79 | b | 83 |
| 14 | d | 85 = 5 . 17 | h | 89 |
| 15 | a | 91 = 7 . 13 | e | 95 = 5 . 19 |
| 16 | g | 97 | b | 101 |
| 17 | d | 103 | h | 107 |
| 18 | a | 109 | e | 113 |
| 19 | g | 115 = 5 . 23 | b | 119 = 7 . 17 |
| 20 | d | 121 = 11 . 11 | h | 125 = 5 . 25 = 5 . 5 . 5 |

| | | | | | |
|---|---|---|---|---|---|
| 21 | a | 127 | | e | 131 |
| 22 | g | 133 = 7 . 19 | | b | 137 |
| 23 | d | 139 | | h | 143 = 11 . 13 |
| 24 | a | 145 = 5 . 29 | | e | 149 |
| 25 | g | 151 | | b | 155 = 5 . 31 |
| 26 | d | 157 | | h | 161 = 7 . 23 |
| 27 | a | 163 | | e | 167 |
| 28 | g | 169 = 13 . 13 | | b | 173 |
| 29 | d | 175 = 7 . 25 = 5 . 35 = 7 . 5 . 5 | | h | 179 |
| 30 | a | 181 | | e | 185 = 5 . 37 |
| 31 | g | 187 = 11 . 17 | | b | 191 |
| 32 | d | 193 | | h | 197 |
| 33 | a | 199 | | e | 203 = 7 . 29 |
| 34 | g | 205 = 5 . 41 | | b | 209 = 11 . 19 |
| 35 | d | 211 | | h | 215 = 5 . 43 |
| 36 | a | 217 = 7 . 31 | | e | 221 = 13 . 17 |
| 37 | g | 223 | | b | 227 |
| 38 | d | 229 | | h | 233 |
| 39 | a | 235 = 5 . 47 | | e | 239 |
| 40 | g | 241 | | b | 245 = 5 . 49 = 7 . 35 = 5 . 7 . 7 |
| 41 | d | 247 = 13 . 19 | | h | 251 |
| 42 | a | 253 = 11 . 23 | | e | 257 |
| 43 | g | 259 = 7 . 37 | | b | 263 |
| 44 | d | 265 = 5 . 53 | | h | 269 |
| 45 | a | 271 | | e | 275 = 5 . 55 = 25 . 11 = 5 . 5 . 11 |
| 46 | g | 277 | | b | 281 |
| 47 | d | 283 | | h | 287 = 7 . 41 |
| 48 | a | 289 = 17 . 17 | | e | 293 |
| 49 | g | 295 = 5 . 59 | | b | 299 = 13 . 23 |
| 50 | d | 301 = 7 . 43 | | h | 305 = 5 . 61 |
| 51 | a | 307 | | e | 311 |
| 52 | g | 313 | | b | 317 |
| 53 | d | 319 = 11 . 29 | | h | 323 = 17 .19 |
| 54 | a | 325 = 25 . 13 = 5 . 5 . 13 | e | 329 = 7 . 47 |
| 55 | g | 331 | b | 335 = 5 . 67 |
| 56 | d | 337 | h | 341 = 11 . 31 |
| 57 | a | 343 = 7 . 49 = 7 . 7 . 7 | e | 347 |
| 58 | g | 349 | b | 353 |
| 59 | d | 355 = 5 . 71 | h | 359 |
| 60 | a | 361 = 19 . 19 | e | 365 = 5 . 73 |
| 61 | g | 367 | b | 371 = 7 . 53 |
| 62 | d | 373 | h | 377 = 13 . 29 |
| 63 | a | 379 | e | 383 |
| 64 | g | 385 = 5 . 77 = 5 . 7 . 11 | b | 389 |
| 65 | d | 391 = 17 . 23 | h | 395 = 5 . 79 |
| 66 | a | 397 | e | 401 |
| 67 | g | 403 = 13 . 31 | b | 407 = 11 . 37 |
| 68 | d | 409 | h | 413 = 7 . 59 |
| 69 | a | 415 = 5 . 83 | e | 419 |
| 70 | g | 421 | b | 425 = 5 . 85 = 17 . 25 = 5 . 5 . 17 |
| 71 | d | 427 | h | 431 |
| 72 | a | 433 | e | 437 = 19 . 23 |
| 73 | g | 439 | b | 443 |
| 74 | d | 445 = 5 . 89 | h | 449 |
| 75 | a | 451 = 11 . 41 | e | 455 = 5 . 91 = 5 . 7 . 13 |
| 76 | g | 457 | b | 461 |
| 77 | d | 463 | h | 467 |
| 78 | a | 469 = 7 . 67 | e | 473 = 11 . 43 |
| 79 | g | 475 = 5 . 95 = 5 . 5 . 19 | b | 479 |
| 80 | d | 481 = 13 . 37 | h | 485 = 5 . 97 |
| 81 | a | 487 | e | 491 |
| 82 | g | 493 = 17 . 29 | b | 497 = 7 . 71 |
| 83 | d | 499 | h | 503 |
| 84 | a | 505 = 5 . 101 | e | 509 |
| 85 | g | 511 = 7 . 73 | b | 515 = 5 . 103 |
| 86 | d | 517 = 11 . 47 | h | 521 |

| | | | | | |
|---|---|---|---|---|---|
| 87  | a | 523            | e | 527 = 17 . 31 |
| 88  | g | 529 = 23 . 23  | b | 533 = 13 . 41 |
| 89  | d | 535 = 5 . 107  | h | 539 = 7 . 77 = 11 . 49 = 7 . 7 . 11 |
| 90  | a | 541            | e | 545 = 5 . 109 |
| 91  | g | 547            | b | 551 = 19 . 29 |
| 92  | d | 553 = 7 . 79   | h | 557 |
| 93  | a | 559 = 13 . 43  | e | 563 |
| 94  | g | 565 = 5 . 113  | b | 569 |
| 95  | d | 571            | h | 575 = 5 . 115 = 5 . 5 . 23 |
| 96  | a | 577            | e | 581 = 7 . 83 |
| 97  | g | 583 = 11 . 53  | b | 587 |
| 98  | d | 589 = 19 . 31  | h | 593 |
| 99  | a | 595 = 5 . 119 = 5 . 7 . 19 | e | 599 |
| 100 | g | 601            | b | 605 = 5 . 121 = 5 . 11 . 11 |

Table VIII. **α** and **β** prime and composite numbers up to **605**.

**8. CONCLUSIONS**

We have described an algorithm to arrange all the positive and negative integer numbers. This array of numbers permits grouping them in six different Classes, **α**, **β**, **γ**, **δ**, **ε**, and **ζ**. Particularly, numbers that belong to Class **α** are defined as **α = 1 + 6 n**, and those of Class **β**, as **β = 5 + 6 n**, where **n = 0, ± 1, ± 2, ± 3, ± 4**, … These two Classes **α** and **β**, contain: i) all prime numbers, except **+ 2**, **–2** and **± 3**, which belong to **ε**, **δ**, and **γ** Classes, respectively, and ii) all the other odd numbers, except those that are multiple of **± 3**, according to the sequence **± 9, ± 15, ± 21, ± 27**, … Besides, products between numbers of the Class **α**, and also those between numbers of the Class **β**, generate numbers belonging to the Class **α**. On the other side, products between numbers of Class **α** with numbers of Class **β** result in numbers of Class **β**. Then, we have demonstrated that both Classes **α** and **β** include: i) all the prime numbers, except **± 2** and **± 3**, and ii) all the products between **α** numbers, as **α. α'**; all the products between **β** numbers, as **β.β'**; and also all the products between numbers of Classes **α** and **β**, as **α.β**, which necessarily are composite numbers, whose factorization is completely determined.

## 9. REFERENCES

[1] "M+D: Otra mirada sobre los números primos" (In spanish), Leopoldo Garavaglia y Mario Garavaglia,
To be presented at the Internacional Conference on Mathematics and Design (M&D 2007), Blumenau, Brazil, 1$^{st}$ to 4$^{th}$ of July, 2007.

## 10. NOTES

Authors can be contacted at: Leopoldo Garavaglia, garavaglia_leo@hotmail.com, and Mario Garavaglia, garavagliam@ciop.unlp.edu.ar